\documentclass{article}\begin{document}\centerline{\bf A Note on Integrals Containing the Univariate Lommel Function}\vskip .7in

\centerline{M.L. Glasser}\vskip .2in

\centerline{Departimento de F\'isica Te\'orica, At\'omica y \"Optica,}
\centerline{Universidad de Valladolid, 47011 Valladolid (Spain)}\vskip .1in

\centerline{Department of Physics}
\centerline{Clarkson University, Potsdam, NY 13699 (USA)}

\vskip 1in
\centerline{\bf ABSTRACT}\vskip .1in
\begin{quote}
This short note investigates a number of index integrals of products of the Lommel functions $s_{\mu,\nu}(a)$ and uncovers an integral relationship. between this function and the Chebyshev polynomials $T_{2n}(x)$.

\end{quote}

\vskip .5in\noindent
Key Words: Lommel Function, Index Integral, Tchebyshev Polynomial, Fourier Inversion.

\vskip .5in\noindent
AMS Classification: 33C10, 44A15, 44A20

\newpage
\centerline{\bf Introduction}\vskip .1in

Index integrals, such as the well-known Kontorovich-Lebedev transform [1]
$$\int_0^{\infty} f(x)K_{ix}(y)dx\eqno(1)$$
 have been investigated and tabulated for over two centuries. A survey of their evaluation and applications can be found in the book
{ \it{ Index Transforms} }[2] by S. Yakubovich.  index integrals of the Lommel functions

$$s_{\mu,\nu}(z)= \frac{z^{\mu+1}}{(\mu+1)^2-\nu^2}\;_1F_2\left(1;\frac{\mu-\nu+3}{2},\frac{\mu+\nu+3}{2};-\frac{z^2}{4}\right) \eqno(2)$$
seem not to  be found in the extensive table by Prudnikov et al.[3] and elsewhere [4] and therefore even a few that can be obtained by unsophisticated means would seem to be worthy of notice, and that is the aim here.

In the next section we present a number of index integrals and sketch their derivation; we follow  this by a few corollary results.

There have been hints in the literature of relations between the Lommel functions and various families of orthogonal polynomials[7]. We verify that this is the case  by investigating a connection between the Lommel function (2) and the Chebyshev polynomials of the first kind and subsequently exploit this nexus to derive a number of sum rules and series, which appear to be new.

\centerline{\bf Results and Derivations}\vskip .1in

\noindent
{\bf Theorem 1.}\vskip .07in
    $$\int_0^{\infty} x\sin(\pi x)s_{-1,x}(a)s_{0,x}(b) dx=\frac{\pi^2}{4}[{\bf H}_0(a-b)-{\bf H}_0(a+b)]\eqno(a)$$

    $$\int_0^{\infty}x \sin(\pi x)s_{-1,x}(a) s_{0,x}(a)dx=-\frac{\pi^2}{4}{\bf H}_0(2a).\eqno(a')$$
    $$\int_0^{\infty}\cos^2(\pi x/2)s_{0,x}(a)s_{0,x}(b)dx=\frac{\pi^2}{8}[J_0(|a-b|)-J_0(a+b)]\eqno(b)$$
    $$\int_0^{\infty}\cos^2(\pi x/2) s_{0,x}^2(a) dx=\frac{\pi^2}{8}(1-J_0(2a)).\eqno(b')$$
    $$\int_0^{\infty}x^2\sin^2(\pi x/2)s_{-1,x}(a)s_{-1.x}(b)=\frac{\pi^2}{8}[J_0(|a-b|)+J_0(a+b)]\eqno(c)$$
    $$\int_0^{\infty} x^2\sin^2(\pi x/2) s_{-1,x}^2(a)dx=\frac{\pi^2}{8}(1+J_0(2a)).\eqno(c')$$
    
    \vskip .1in\noindent
    {\bf Proof}: These are all consequences of the Fourier transforms[ 5,(1.7.(49),(50))] whose inversion gives
    
    $$\sin(a\cos(x))\theta\left(\frac{\pi}{2}-x\right)=\frac{2}{\pi}\int_0^{\infty}\cos(xy)\cos\left(\frac{\pi y}{2}\right)s_{0,y}(a)dy\eqno(3a)$$
    $$\cos(b\cos x))\theta\left(\frac{\pi}{2}-x\right)=-\frac{2}{\pi}\int_0^{\infty}\cos(xy)y\sin\left(\frac{\pi y}{2}\right)s_{-1,y}(b)dy\eqno(3b)$$
    where $\theta$ denotes the unit step function: 1 for positive argument, 0 for negative. By Parseval's identity  and the double angle identity for the sine we have
    $$\int_0^{\pi/2}\sin(a\cos x)\cos(b\cos x)dx=-\frac{1}{\pi}\int_0^{\infty} x\sin(\pi x)s_{-1,x}(a)s_{0,x}(b) dx.\eqno(4)$$
    Since
    
    $$-\pi\int_0^{\pi/2}\sin(a\cos x)\cos(b\cos x)dx
    =\frac{\pi^2}{4}[{\bf H}_0(a+b)+{\bf H}_0(a-b)]\eqno(5)$$
    one has (a) and (a'). The derivations of the pairs $b$ and $c$ proceed similarly.

    Next, the Fourier inversions of (3a,b) with $x=cos^{-1}u$ read
    $$\int_0^1\sin(au)\frac{\cos(y\cos^{-1}u)}{\sqrt{1-u^2}}du=\cos\left(\frac{\pi y}{2}\right)s_{0,y}(a)\eqno(6a)$$
    $$\int_0^1\cos(bu)\frac{\cos(y\cos^{-1}u)}{\sqrt{1-u^2}}du=-y\sin\left(\frac{\pi y}{2}\right)s_{-1,y}(b).\eqno(6b)$$
    so, again by Fourier inversion,
    $$\frac{\cos(y\cos^{-1} u)}{\sqrt{1-u^2}}\theta(1-u)=\frac{2}{\pi}\cos\left(\frac{\pi y}{2}\right)\int_0^{\infty}\sin(ut)s_{0,y}(t)dt ,\quad u\ne0\eqno(7a)$$
    $$\frac{\cos(y\cos^{-1} u)}{\sqrt{1-u^2}}\theta(1-u)=-\frac{2}{\pi}y\sin\left(\frac{\pi y}{2}\right) \int_0^{\infty}\cos(ut)s_{-1,y}(t)dt.\eqno(7b)$$
    Therefore, with $y=2n, 2n+1$, respectively, one has the explicit connection with Tchebyshev polynomials
    
    \noindent
    {\bf Theorem 2}
    $$\frac{T_{2n}(u)}{\sqrt{1-u^2}}\theta(1-u)=(-1)^n\frac{2}{\pi}\int_0^{\infty}\sin(ut)s_{0,2n}(t)dt\eqno(8a)$$
    $$\frac{T_{2n+1}(u)}{\sqrt{1-u^2}}\theta(1-u)=(-1)^{n+1}(2n+1)\frac{2}{\pi}\int_0^{\infty}\cos(ut)s_{-1,2n+1}(t)dt\eqno(8b)$$
  
    $$s_{0,2n}(t)=(-1)^n\int_0^1\sin(ut)\frac{T_{2n}(u)}{\sqrt{1-u^2}}du\eqno(9a)$$
    $$s_{-1,2n+1}(t)=\frac{(-1)^{n+1}}{2n+1} \int_0^1\cos(ut)\frac{T_{2n+1}(u)}{\sqrt{1-u^2}}du.\eqno(9b)$$
    
    \vskip .2in
    
   \centerline{\bf Discussion}\vskip .1in
   
   The results of Theorem 1 can be extended to other integer values of the first index by means of documented[6] recursion relations for the Lommel function. For example
   $$s_{0,x}(a)=\frac{a-s_{2,x}(a)}{1-x^2},\quad  s_{-1,x}(a)=x^{-2} s_{1,x}(a)\eqno(10)$$
   $$s_{m,x}(a)=\frac{a}{2x}[(m+x-1)s_{m-1,x-1}(a)-(m-x=1)s_{m-1,x+1}(a)]\eqno(111)$$
   $$\frac{d}{da}s_{m,x}(a)=\frac{1}{2}[(m+x-1)s_{m-1,x-1}(a)+(m-x-1)s_{m-1,x+1}(a)].\eqno(12)$$

   Since the Tchebyshev functions are related to other families of orthogonal polynomials[7] it is possible to connect these with the Lommel function as well. For example, 
   $$U_{2n}(x)=\frac{(-1)^n}{\sqrt{1-x^2}}T_{2n+1}(\sqrt{1-x^2})\eqno(13)$$
   gives for the Tchebyshev polynomial of the second kind
   $$U_{2n}(x)=-(2n+1)\int_0^{\infty}\cos(t\sqrt{1-x^2})s_{-1,2n+1}(t)dt\eqno(14)$$
   Finally, one can see from (9a) and (9b) that the  differentialtion formula
   $$2s'_{\mu,\nu}(x)=(\nu-1)s_{\mu-1,\nu-1}(x)-(\nu+1)s_{\mu-1,\nu+1}(x)\eqno(15a)$$
    is equivalent to the Chebyshev identity
    $$T_{2n+1}(x)+T_{2n-1}(x)=2x T_{2n}(x)\eqno(15b)$$
    for $\mu=0$, $\nu=2n$.There is an abundance of the former in the compilation [8], for example
    
    It thus appears that almost every relation satisfied by the Tchebyshev polynomials is mirrored by one for $s_{\mu,\nu}(t)$.   For example, from[8,5.7.1(3)], (8a) and Fourier inversion, one has
    $$\sum_{k=0}^n (-1)^k s_{0,2k}(x)=\frac{\pi}{4}{\bf H}_0(x)+\frac{1}{2}\int_0^1\sin(xt)\frac{U_{2n}(t)}{\sqrt{1-t^2}}dt,\eqno(16)$$
    where $U_{2n}$ denotes the Tchebyshev polynomial of the second kind. Similarly, from [8,6.12.2(1)]
    $$\sum_{k=0}^{\infty}J_{2k}(w)s_{0,2k}(t)=  \frac{\pi}{4}  [ {\bf H}_0  ( t-w.)+{\bf H}_0(t+w)+2J_0(w) {\bf H}_0(t)].\eqno(17)$$
   
   \vskip .2in
     \centerline{\bf Acknowlegements}
   I wish to that Prof. S. Yakubovich,  Prof, Roy Hughes  and Dr. Michael Milgram for valuable correspondence.
I also appreciate the support of two research projects: QCAYLE, funded by the European Union, NextGenerationEU, and PID2020-113406GB-I0 funded by the MCIN of Spain. 
   
   \newpage

   \centerline{\bf References}\vskip .1in
   \noindent
   [1] A.\'Erdelyi et al., {\it Tables of Integral Transforms, Vol.2}  [McGraw-Hill, NY.(1952); Chap.XII]
   
  \noindent
  [2]  S. B. Yakubovich, {\it Index Transforms}, [World Scientific Publishers, Singapore 1996]
  
  \noindent
  [3] A. P. Prudnikov et al., {\it Table of Integrals, Series and Products Vol3},{Nauka Publishers, Moscow (1983); Sec.2.9].
  
  \noindent
  [4] S.B. Yakubovich, Integral Transformations by the Index of Lommel's Function, Periodica Mathematica Hungarica {\bf 16}, 223-233(2003).
  
  \noindent
  [5] A. \'Erdelyi et al.{\it Tables of Integral Transforms, Vol1.},[McGraw-Hill, NY (1953)]
  
  \noindent
  [6] A. \'Erdelyi, et al.{ \it Higher Transcendental Functions, Vol.2}[McGraw-Hill, 1953; Sec. 7.3.5]
  
  \noindent
  [7] H.T. Kpelink and W.A. van Assche, Orthogonal Polynomials and Laurent Polynomials Related to the Hahn-Exton q-Besseel Function. arXiv:math 9502227 v1(1995)
  
  \noindent
  [8]  Yuri A. Brychkov, {\it Handbook of Special Functions},[CRC Press, Chapman and Hall Pub.Boca Raton, Fla.(USA) 2008]

\end{document}